\documentclass[12pt]{amsart}
\usepackage{amsmath,amscd,amssymb,amsthm,array}
\usepackage{cancel}
\usepackage{mathtools}

\usepackage{stmaryrd}

\newcommand{\llb}{\llbracket}
\newcommand{\rrb}{\rrbracket}

\usepackage{hyperref}
\usepackage{tikz-cd}

\newcommand {\del}{{\partial}}

\usepackage{DLdef1}

\let\ssec\subsection

\renewcommand {\ssbegin}[2][*]
 {\refstepcounter{subsection}%
\if#1*
\addcontentsline{toc}{subsection}{\thesubsection.\hskip 1pc #2}%
\else
\addcontentsline{toc}{subsection}{\thesubsection.\hskip 1pc #2. #1}%
\fi
 \def \secno {\gdef \secno {}{\ssecfont
\thesubsection.\hskip 2ex}%
 }%
 \begin{#2}}

\renewcommand {\sssbegin}[2][*]
 {\refstepcounter{subsubsection}
\if#1*
\addcontentsline{toc}{subsubsection}{\thesubsubsection.\hskip 1pc #2}%
\else
\addcontentsline{toc}{subsubsection}{\thesubsubsection.\hskip 1pc #2. #1}
\fi
 \def \secno {\gdef \secno {}{\ssecfont \thesubsubsection.\hskip 2ex}%
 }%
 \begin{#2}}

\renewcommand {\parbegin}[2][*]
 {\refstepcounter{paragraph}
\if#1*
\addcontentsline{toc}{paragraph}{\theparagraph.\hskip 1pc #2}%
\else
\addcontentsline{toc}{paragraph}{\theparagraph.\hskip 1pc #2. #1}
\fi
 \def \secno {\gdef \secno {}{\ssecfont \theparagraph.\hskip 2ex}%
 }%
 \begin{#2}}

\title{Gordan-Rankin-Cohen operators on the spaces of weighted densities in superdimension $1\vert 1$}

\author[V.~Bovdi, D.~Leites]{Victor Bovdi${}^{a,*}$, Dimitry Leites${}^{b}$}
\address{
 ${}^*$Corresponding author\\
 ${}^a$Department of Mathematics, UAEU, Al Ain, UAE; vbovdi@gmail.com\\
${}^b$Department of Mathematics Stockholm University, Albanov\"agen 28, SE-114 19, Stockholm,
 Sweden; dimleites@gmail.com\\}

\begin{document}

\begin{abstract} The modular forms and weighted densities over the 1-dimensional manifold $M$ are transformed ``alike" under the group 
of linear fractional changes of coordinates, so the classifications of differential operators between spaces of (A) modular forms and (B) weighted densities were identified in various works. Actually, these problems are different. Here, we solve  problem~ B for superstrings in superdimension $(1\vert 1)$ --- superizations of the result of arXiv:2404.18222. 
Open problems are offered.
\end{abstract}

\keywords{Lie superalgebra, invariant differential operator,
Gordan transvectant, Rankin-Cohen bracket, weighted density, superstring}

\makeatletter
\@namedef{subjclassname@2020}{\textup{2020} Mathematics Subject Classification}
\makeatother
\subjclass[2020]{Primary 17B10 Secondary 53B99, 32Wxx}

\maketitle

\markboth{{\itshape V.~ Bovdi\textup{,} D. ~Leites}}
{{\itshape Gordan-Rankin-Cohen operators on superstrings of superdimension $1\vert 1$}}

\section{Introduction} \label{Intro}

\subsection{Weighted densities vs. modular forms} Recall that on any finite-dimensional supermanifold $\cM$, the elements of the space $\cF_w$ of weighted densities of \textit{weight} $w\in\Cee$ are expressions of the form $f(\vvol)^w$, where $f\in\cF:=\cF_0$ is a~function  and $\vvol$ is a~volume element. The Lie derivative $L_X$ along the vector field $X$ with divergence $\Div X$ (for its precise form, see eq.~\eqref{Div}) acts by the formula
\be\label{vvol}
L_X(f(\vvol)^w)=(X(f)+wf\Div X)(\vvol)^w.
\ee
Equivalently, the weighted densities can be described in terms of any subsupergroup $G\subset \Diff (\cM)$ of diffeomorphisms $A\in G$ acting on $\cM$ by changes of coordinates $z\to \widetilde z:=A(z)$:
\be\label{vvol1}
A(f(z)(\vvol(z))^w)\tto f(\widetilde z)(\vvol(\widetilde z))^w=f(\widetilde z)\left(\ber\left(\nfrac{\partial\tilde z}{\partial z}\right)\vvol(z)\right)^w,
\ee
where $\ber$ is the berezinian (superdeterminant) 
of the Jacobi matrix of the change of coordinates.

On the line, $\vvol$ is $dz$. Since $d\tilde z=\nfrac{\partial\tilde z}{\partial z} dz$ and
\be\label{1}
\text{$\nfrac{\partial\tilde z}{\partial z}=(cz+d)^{-2}$ \quad for $\tilde z=\nfrac{az+b}{cz+d}$\quad and \quad $\begin{pmatrix}a&b\\
c&d\end{pmatrix}\in \SL(2;\Cee)$,}
\ee
then, according to definition \eqref{vvol1}, and since on manifolds $\ber=\det$, we replace $z\mapsto \tilde z$ in $f(z)\vvol(z)^{-k/2}$ and get
\be\label{confuse} 
f(\widetilde z)(\vvol(\widetilde z))^{-k/2}=f(\widetilde z)\left(\ber\left(\nfrac{\partial\tilde z}{\partial z}\right)\vvol(z)\right)^{-k/2}=f\left (\nfrac{az+b}{cz+d}\right)(cz+d)^{k}\vvol(z)^{-k/2}.
\ee

Recall that a~holomorphic function $\varphi$ on the upper half-plane $\cH$ is called a~\textit{modular form of weight}~ $k$ (``typically, a~positive integer", see \cite{W}) if under the $\SL(2;\Zee)$-action on $\cH$ by  linear fractional changes of coordinates the function $\varphi$ is transformed as follows: 
\be\label{lf}
\varphi\left(\nfrac{az+b}{cz+d}\right)=\varphi(z)(cz+d)^{k}\qquad \text{for any}\quad \begin{pmatrix}a&b\\
c&d\end{pmatrix}\in\SL(2;\Zee),
\ee
and, moreover, $\varphi$ satisfies the following growth conditions, where $z=x+iy$, see \cite[p.\,5]{Z2}
\be\label{grCond}
\text{$\varphi(x + iy) \sim O(1)$ as $y \to\infty$ and $\varphi(x + iy) \sim O(y^{-k})$ as $y \to 0$.}
\ee

Note that
\begin{equation}\label{wden}
\begin{minipage}[c]{13cm}
If $\varphi$ in $\varphi\,\vvol^{k/2}$ satisfies \eqref{lf},
then $\varphi\,\vvol^{k/2}$
is $\SL(2;\Zee)$-invariant.\\
\end{minipage}
\end{equation}

\sssec{Vital difference and a cause of confusions}\label{DiffConf} The modular forms and weighted densities over the 1-dimensional manifold $M$ are transformed ``alike" under the group 
of linear fractional changes of coordinates, so the following two problems are sometimes identified 
\begin{equation}\label{2prob}
\begin{minipage}[c]{13cm}
``classify differential operators between spaces of\\
(A) modular forms,\\ (B) weighted densities".
\end{minipage}
\end{equation}
I.~Shchepochkina observed that under the $\SL(2; \Zee)$-action the modular form is \textbf{not} ``transformed 
as the coefficient of the volume element in the weighted density", as one might find in, e.g.,~\cite{OR, GO, BLO, BoLe} and even in \cite[p.~41, lines 8--10]{CMZ} corrected in \cite{LSh}. Indeed, 
\begin{equation}\label{2probl}
\begin{minipage}[c]{13cm}
the modular form $\varphi$ \textbf{itself} is multiplied by $(cz + d)^k$, see eq.~\eqref{lf}, whereas one has to additionally perform a~change of variables in the coefficient $f$ of the volume element, see eq.~\eqref{confuse}.
\end{minipage}
\end{equation} 

Besides, the spaces $\cF_w$ of any weight $w\in\Cee$ are $\cF$-modules of rank 1, while thanks to additional condition \eqref{grCond}, the spaces $M_k$ of modular forms of weight $k\in\Zee_+$ are finite-dimensional over $\Cee$ (see \cite[eq.(7)]{Z2}):
\be\label{dim}
\dim M_k=
\begin{cases}0, &\text{\quad for $k < 0$ or $k$ odd};\\
\lfloor k/12\rfloor + 1, &\text{\quad if $k \not\equiv 2\pmod{12}$};\\
\lfloor k/12\rfloor, &\text{\quad if $k \equiv 2 \pmod{12}$}.
\end{cases}
\ee

\textbf{The unary  $\SL(2;\Zee)$-invariant differential operators} between spaces of modular forms, as well as the $\SL(2;\Cee)$-invariant differential operators between spaces of weighted densities, are called \textit{Bol} operators, see \cite{Bol, Gi}. In \cite{BLS2}, it was shown that the description of Bol operators in the sense of problem B, see eq.~\eqref{2prob}, on supermanifolds without any structure preserved is a~wild problem, except for superstrings of superdimension $(1|1)$ in which case the answer was given. For the superstrings of any dimension $(1|n)$ with a~contact structure, the answer is given in \cite{BLS1}; some of these super Bol operators are very interesting.

For the \textbf{unary} (Bol) operators  on 1-dimensional and $(1|1)$-dimensional superstrings with a contact structure,  despite the difference in the dimensions and nature of the modular forms and weighted densities, 
solutions to problem A are in one-to-one correspondence with solutions to problem B: compare the answers in~\cite{Bol,Gi} with those in \cite{BLS1}. 

\textbf{The bilinear $\SL(2;\Zee)$-invariant differential operators} between spaces of modular forms (solutions to problem A, see eq.~\eqref{2prob}) are called \textit{Gordan transvectors} or \textit{Rankin--Cohen} brackets, were introduced and rediscovered in different settings in \cite{Gor, R, Co, Z, JP}. 

Let us call any solution to problem~B, see eq.~\eqref{2prob}, understood as any bilinear operator invariant under any (not necessarily simple and maximal) subalgebra $\fm$ of the Lie (super)algebra $\fg(\cM)$ of vector fields  on the (super)manifold $\cM$ preserving a structure on $\cM$  a~ GRC \textbf{operator} to distinguish such operators from the GRC \textbf{brackets}.

For the \textbf{binary} operators, to every solution to problem A a~solution to problem B corresponds, but the latter are more abundant: if $\dim\cM= 1$, compare the brackets in~Theorem~in \cite[p.\,57]{Z2}  and the brackets considered in \cite{OR} on the one hand with the brackets listed in \cite{BoLe} on the other hand; in dimension $(1|1)$, compare lists of the brackets in \cite{Gi, CMZ, GO} with our Theorem~\ref{ThGRC2}.

\sssec{On multidimensional generalizations and superizations} \textbf{On 1-dimensional manifolds}, the problems A and B are related sharing similar transformation rules \eqref{lf} and  \eqref{confuse}, but still differ: see \eqref{2probl}. 

\textbf{Problem B}. \textbf{On superstrings of any dimension}, the unary (resp., binary) differential operators between spaces of weighted densities invariant under the maximal and simple subalgebra of all vector fields $\fvect(1|N)$ or under the maximal and simple subalgebra of  contact vecotr fields $\fk(1|N)$ were analogously called
\textit{Bol} (resp., \textit{GRC}) operators, see \cite{Gi, GO, BLS1, BLS2} and references therein. Let $\fg:=\fg(\cM)=\oplus\ \fg_i$ be a~$\Zee$-graded Lie superalgebra of polynomial (or formal) vector fields on the $(1|n)$-dimensional superstring $\cM$ with even coordinate $t$ and odd coordinates $\theta:=(\theta_1, \dots,  \theta_n)$. Let $V$ be a~$\fg_0$-module; by setting $\fg_iV=0$ for $i>0$ we turn $V$ into a~module over $\fg_{\geq0}:=\oplus_{i\geq0}\ \fg_i$, and hence over $U(\fg_{\geq0})$. The \textit{space of tensor fields of type} $V$ is defined to be
\be\label{T(V)}
T(V):=\Hom_{U(\fg_{\geq0})}(U(\fg), V)\simeq\Cee[[t, \theta]]\otimes V.
\ee

If the non-trivial irreducible $\fg_0$-module $V$ is 1-dimensional, the elements of $T(V)$ are analogs of weighted densities. For the three Lie superalgebras $\fg$, namely: $\fvect(1|1)$ of all polynomial vector fields, the Lie superalgebra  $\fk(1|1)$ of fields preserving a distribution given by the contact form with even "time", and the Lie superalgebra $\fm(1)$ of fields preserving a distribution given by the contact form with an odd  ``time" (for details, see \cite{LSh}), all non-trivial irreducible $\fg_0$-modules~ $V$ are 1-dimensional. In this article, we have solved problem B in the first two of these three cases.

\section{The GRC operators on the $n=1$ superstring with a contact structure} \label{SS11c}

Note that in presence of a contact structure, it is natural to express the volume element $\vvol$ as a power of the contact form $\alpha_1$. 

\ssec{$\cS_c^{1|n}$ and $\alpha_{1}$}\label{SSalpha1} Let a~distribution on the $(1|n)$-dimensional superstring $\cS_c^{1|n}$ with an even coordinate $t$  and odd coordinates $\theta:=(\theta_1, \dots,  \theta_n)$ be given by the 1-form $\alpha_1:=dt+\sum \theta_i d\theta_i$.
Let $\cF$ be the supercommutative superalgebra of functions (polynomials or Laurent polynomials) on $\cS_c^{1|n}$ and $\cF_\mu:=\cF\alpha_1^{\mu/2}$ for $n>0$ the space of $\mu$-densities for any $\mu\in\Cee$. Let $E:=\sum \theta_{i}\partial_{\theta_i}$. The Lie superalgebra $\fk(1|n)$ is spanned by the contact vector fields  
\be\label{Kf}
K_{f}:=\textstyle(2-E)(f)\partial_t- \displaystyle (-1)^{p(f)}\sum \left(\pderf{f}{\theta_i}
 - \theta_i\pderf{f}{t}\right) \partial_{\theta_i} \qquad \text{for any}\quad f\in\cF
\ee
that preserve the distribution given by $\alpha_1$: indeed, $L_{K_{f}}(\alpha_1)=2\del_t(f)\alpha_1$.
The contact bracket of generating functions is given by the formula
\be\label{kb}
\{f,\, g\}_{k.b.}:=\textstyle(2-E)(f)\pderf{g}{t}-\pderf{f}{t}(2-E)(g)-
\displaystyle (-1)^{p(f)}\sum \pderf{f}{\theta_i} \pderf{g}{\theta_i}\qquad \text{for any}\quad f, g\in\cF.
\ee

\ssec{Classification of $\fosp(1|2)$-invariant GRC operators} Let $\fg:=\fk(1|1)$ be generated by polynomial functions. Set $\deg \theta=1$, $\deg t=2$, hence $\deg K_f=\deg f-2$ and $\fg_0=\Cee K_t$. The embedded subalgebra $\fosp(1|2)\subset\fk(1|1)$ is generated by $\nabla_-:=K_{\theta}=\theta\del_t+\del_{\theta}$ and $\nabla_+:=K_{t\theta}$. The eigenvector with respect to $K_t$ will be called a~\textit{weight} vector. 

Let $V$ be a~$\fg_0$-module; setting $\fg_iV=0$ for $i>0$ we turn $V$ into a~module over $\fg_{\geq0}:=\oplus_{i\geq0}\ \fg_i$. 
We define the space dual to $T(V)$, see eq.~\eqref{T(V)}, by setting
\[
I(V^*):=U(\fg)\otimes_{U(\fg_{\geq0})} V^*\simeq\Cee[K_1,K_{\theta}]\otimes V^*\sim \Cee[K_{\theta}]\otimes V^*.
\]
Note that, strictly speaking, the last $\sim$ instead of $\simeq$ above and in eqs.~\eqref{hmm}, \eqref{sim} reflect a~slight ``rudeness of speech", since, speaking correctly, we mean an isomorphism with the quotient space
\[
\Cee[K_1, K_\theta]/(K_1-2(K_\theta)^2)\simeq \Cee[K_\theta]\qquad \text{since} \quad K_1=[K_\theta, K_\theta]=2(K_\theta)^2.
\]

Since $\fg_0=\fgl(1)$, all finite-dimensional irreducible $\fg_0$-modules are
1-dimensional; let $V^*_1$ and $V^*_2$ be such
modules; let them be even. Let $v\in V^*_1$ and $w\in V^*_2$ be nonzero vectors of
weight $\mu_1$ and $\mu_2$, respectively, i.e.,
\[
K_t(v)=\mu_1v, \qquad K_t(w)=\mu_2w.
\]

Setting $X(V^*_1)=X(V^*_2)=0$ for any $X\in\fg_{>0}$, where $\fg_{> 0}:=\oplus_{i> 0}\ \fg_i$, we turn $V^*_1$ and $V^*_2$ into modules over $\fg_{\geq 0}:=\oplus_{i\geq 0}\ \fg_i$. Let $K_f'$ and $K_f''$ be copies of $K_f$. Set 
\be\label{hmm}
\begin{split}
I(V^*_1, V^*_2):&=(U(\fg)\otimes_{U(\fg_{\geq 0})}V^*_1)\otimes (U(\fg)\otimes_{U(\fg_{\geq 0})}V^*_2)\\
&\simeq (\Cee[K_1', K_\theta']\otimes V^*_1)\bigotimes (\Cee[K_1'', K_\theta'']\otimes V^*_2)\\
&{\sim} (\Cee[K_\theta']\otimes V^*_1)\bigotimes (\Cee[K_\theta'']\otimes V^*_2).
\end{split}
\ee
Note that $\fg_>$ is generated by $K_{t\theta}$ and $K_{t^2\theta}$. Following Rudakov (see \cite{Ru}) recall that a~
\textit{singular vector} in ${I(V_1)\otimes I(V_2)}$  or $I(V)$ is a~weight vector relative $K_t$ that vanishes under ${\fg_{>0}:=\oplus_{i>0}\, \fg_i}$.
Let us describe the $\fosp(1|2)$-\textit{singular} vectors\footnote{Note that in \cite{BoLe}  the term \textit{singular} vector was carelessly applied to both singular vectors and to what should have been called $\fpgl(2)$-\textit{singular} vectors.}
$v_k\in I(V^*_1, V^*_2)$, i.e., the weight vectors relative $K_t$ killed only by $\nabla_+=K_{t\theta}$.
We express these $\fosp(1|2)$-singular vectors as sums of summands with certain coefficients selected for convenience of calculations, as in \cite{Gr}:
\[
\begin{split}
v_{2n}=\textstyle\mathop{\sum}\limits_{0\leq i\leq n}\nfrac{1}{i!(n-i)!}c_i(K_1')^iv&\otimes(K_1'')^{n-i}w\\
&+\textstyle\mathop{\sum}\limits_{0\leq i\leq n-1}e_iK_\theta'(K_1')^iv\otimes K_\theta''(K_1'')^{n-i-1}w;\\
v_{2n+1}=\textstyle\mathop{\sum}\limits_{0\leq i\leq n}a_iK_\theta'(K_1')^iv&\otimes(K_1'')^{n-i}w\\
&+\textstyle\mathop{\sum}\limits_{0\leq i\leq n}\nfrac{1}{i!(n-i)!}b_i(K_1')^iv\otimes K_\theta''(K_1'')^{n-i}w.\\
\end{split}
\]
Since 
\be\label{volISdtheta}
\text{on $\cS_c^{1|n}$ we have $\vvol\simeq [d\theta]\pmod {\cF\alpha}$,}
\ee 
it follows that 
\be\label{sim}
(\cF_{-\mu_1})^*\simeq I(V^*_1)\sim \Cee[K_\theta]\otimes V^*_1
\ee 
implying an interpretation of $\fosp(1|2)$-singular vectors $v_n$ in terms of GRC operators.

Note that $\{t\theta, 1\}=-2\theta$, $\{\theta, \theta\}=1$ and $\{t\theta, \theta\}=t$, see formula~\eqref{kb}. Then,  
\[
\begin{split}
\nabla_+((K_1')^iv\otimes (K_1'')^{n-i}w)=&-2iK_\theta' (K_1')^{i-1}v\otimes (K_1'')^{n-i}w\\
&-2(n-i)(K_1')^iv\otimes K_\theta'' (K_1'')^{n-i-1}w;\\
\end{split}
\]
\[
\begin{split}
\nabla_+(K_\theta' (K_1')^iv\otimes& K_\theta'' (K_1'')^{n-i-1}w)\\
=&(\mu_1-i)(K_1')^iv\otimes K_\theta'' (K_1'')^{n-i-1}w\\
&+ (-\mu_2+3(n-i-1)) K_\theta' (K_1')^iv\otimes (K_1'')^{n-i-1}w;\\
&\\
\nabla_+(K_\theta' (K_1')^iv\otimes& (K_1'')^{n-i}w)\\
=&(\mu_1-i)(K_1')^iv\otimes (K_1'')^{n-i}w+ 2(n-i)K_\theta'(K_1')^iv\otimes K_\theta'' (K_1'')^{n-i-1}w;\\
&\\
\nabla_+((K_1')^iv\otimes& K_\theta'' (K_1'')^{n-i}w)\\
=&-2iK_\theta' (K_1')^{i-1}v\otimes K_\theta'' (K_1'')^{n-i}w+ (\mu_2-(n-i))(K_1')^iv\otimes (K_1'')^{n-i}w.\\
\end{split}
\]
The summands involving $(K_1')^{-1}$ and $(K_1'')^{-1}$ appear only with zero coefficients, so these summands make sense. We consider the two cases.

\underline{Case (i) $\nabla_+(v_{2n})=0$}. The condition
\[
\begin{split}
0&=\nabla_+v_{2n}\\
&=\textstyle \mathop{\sum}\limits_{0\leq i\leq n}\nfrac{c_i}{i!(n-i)!}\left(-2iK_\theta' (K_1')^{i-1}v\otimes (K_1'')^{n-i}w-2(n-i)(K_1')^iv\otimes K_\theta'' (K_1'')^{n-i-1}w\right)\\
&\quad +\textstyle \mathop{\sum}\limits_{0\leq i\leq n-1}\nfrac{e_i}{i!(n-1-i)!}((\mu_1-i)(K_1')^iv\otimes K_\theta'' (K_1'')^{n-i-1}w\\
&\qquad\qquad + (-\mu_2+3(n+i+1)) K_\theta' (K_1')^iv\otimes (K_1'')^{n-i-1}w)
\end{split}
\]
is equivalent to the following system
\be\label{Sect1}
\begin{split}
\textstyle\nfrac{(\mu_1-i)}{i!(n-1-i)!}e_i&=\textstyle\nfrac{2(n-i)}{i!(n-i)!}c_{i};\\
\textstyle\nfrac{(3(n+i+1) -\mu_2)}{i!(n-1-i)!}e_i&=\textstyle\nfrac{2(i+1)}{(i+1)!(n-i-1)!}c_{i+1}, \qquad \text{where}\quad i=0,\dots, n-1.
 \end{split}
\ee
Hence,
\be\label{Sect2}
\big(3(n+i+1) -\mu_2\big)e_{i}+\big((i+1)-\mu_1\big)e_{i+1}=0, \qquad \text{where}\quad i=0,\dots, n-1.
\ee
We denote by $M$ the $n\times (n+1)$ matrix whose elements are the coefficients of system \eqref{Sect2}. Clearly, $M$ can have non-zero elements only on the main diagonal $D=
(A_0,\ldots, A_{n-1})$, where $A_i=3(n+i+1) -\mu_2$,
and on the diagonal just above the main one $(B_0,\ldots, B_{n-1})$, where $B_i=(i+1)-\mu_1$. Obviously, at most one $A_i$ (resp., $B_i$) can vanish. Indeed, if $A_i=A_j=0$ for some $i\not=j$, then $i=j$, which is a~contradiction. Moreover, by condition~ \eqref{Sect2}, if $A_i=0$, then $B_i=0$, and vice versa.

\underline{Case (i)1: $D$ has no zero components}. Then,  from \eqref{Sect2} we have
\be\label{Hek00}
\begin{split}
e_{i}&=(-1)^{n-i}\textstyle\prod_{k=i}^{n-1}\nfrac{B_k}{A_k}e_{n}=\prod_{k=i}^{n-1}\nfrac{\mu_1-(k+1)}{3(n+k+1) -\mu_2}e_{n};\\
c_{i}&=\textstyle\nfrac{\mu_1-i}{2}\prod_{k=i}^{n-1}\nfrac{\mu_1-(k+1)}{3(n+k+1) -\mu_2}e_{n}, \qquad \text{where}\quad i=0,\ldots, n-1.
\end{split}
\ee
Hence, $\dim(\Ker(\nabla_+))=1$ and each $e_i$ and $c_i$ depends on $e_{n}$; so, up to a~factor, $v_i$ depends on~ $e_{n}$.

\underline{Case (i)2: $D$ has a~zero component}. Let $A_j=B_j=0$ for some $j\in\{ 0, \dots, n-1\}$.
Clearly, $\mu_1=j+1$ and $\mu_2=3(n+j+1)$, so $A_i=3(i-j)\not=0$ and $B_i=i-j\not=0$ for all $i\in\{ 0, \dots, n-1\}\setminus \{j\}$. The system \eqref{Sect2} is equivalent to the following system:
\[
\begin{cases}
3e_{i}+ e_{i+1}&=0, \qquad \text{where}\quad i\in\{ 0,\dots, j-1\};\\
0e_{j}+0e_{j+1}&=0;\\
3e_{i}+ e_{i+1}&=0, \qquad \text{where}\quad i\in\{ j+1,\dots, n-1\},\\
\end{cases}
\]
whose solutions depend on two independent parameters $e_j$ and $e_{n}$. Consequently,
\be\label{FinalEquations_01}
\begin{cases}
\begin{split}
e_{i}&=(-3)^{i-j}e_{j};\quad c_{i}=\textstyle\nfrac{j-i+1}{2}(-3)^{i-j}e_{j}, \qquad \text{where} \quad i=0,\ldots, j-1;\\
e_{l}&=(-{3})^{l-n}e_{n};\quad
c_{l}=\textstyle\nfrac{n-l+1}{2}(-{3})^{l-n}e_{n}, \qquad \text{where} \quad l=j+1,\ldots, n-1.
\end{split}
\end{cases}
\ee
Note that if $j=0$, then we ignore the top line of system \eqref{FinalEquations_01} whereas if $j=n-1$, then we ignore the bottom line of system \eqref{FinalEquations_01}.

\underline{Case (ii) $\nabla_+(v_{2n+1})=0$}. From
\[
\begin{split}
0=&\nabla_+v_{2n+1}\\
=&\textstyle \mathop{\sum}\limits_{0\leq i\leq n}{a_i}\left((\mu_1-i)(K_1')^iv\otimes (K_1'')^{n-i}w+ 2(n-i)K_\theta'(K_1')^iv\otimes K_\theta'' (K_1'')^{n-i-1}w\right)\\
&+\textstyle \mathop{\sum}\limits_{0\leq i\leq n-1}{b_i}\left(-2iK_\theta' (K_1')^{i-1}v\otimes K_\theta'' (K_1'')^{n-i}w+ (\mu_2-n+i)(K_1')^iv\otimes (K_1'')^{n-i}w\right )
\end{split}
\]
it follows that
\be\label{Sect4}
\begin{split}
(i-\mu_1)a_i&=(\mu_2-n+i)b_i, \qquad\qquad \text{where}\quad i=0,\dots, n-1;\\
a_i&=\textstyle\nfrac{i+1}{n-i}b_{i+1}, \quad \qquad\qquad \quad \text{where}\quad i=0,\dots, n-2;\\
a_{n-1}&=a_n(\mu_1-n)=0.\\
 \end{split}
\ee

Consider the following subcases:

\underline{Case (ii)1. Let $\mu_2-n+i\not=0$ for all $i\in\{0,\dots, n-1\}$.} Then,  for $i=n-1$ from the 1st equation of \eqref{Sect4} we deduce that $(\mu_2-1)b_{n-1}=0$ and $a_0=\cdots=a_{n-1}=0$ by the 1st and 2nd equations of \eqref{Sect4}. It follows that $b_0=\cdots=b_{n-1}=0$ by the 1st equation of \eqref{Sect4} and $a_n(\mu_1-n)=0$. Consequently,
\be\label{Case1}
\text{if $\mu_2\not\in\{1,\dots, n\}$ and $\mu_1=n$, then there is a~unique (up to a~factor) solution.}
\ee

\underline{Case (ii)2. Let $\mu_2-n+k=0$ for some $k\in\{0,\dots, n-1\}$ and $\mu_1\not\in \{0,\dots, k-1\}$.} Similarly to the previous case, we obtain that
\[
\begin{split}
a_{k}=a_{k+1}&=\cdots=a_{n-1}=(\mu_1-n)a_n=0;\\
b_{k+1}&=b_{k+2}=\cdots=b_{n-1}=0.
\end{split}
\]
Clearly, if $b_k\not=0$, then $a_s\not=0$ for all $s\in\{0,\ldots, k-1\}$ by the 2nd equation of \eqref{Sect4}, so we have the following solution:
\be\label{Case2}
\begin{split}
&b_i=\textstyle\nfrac{(i-\mu_1)(i+1)}{(\mu_2-n+i)(n-i)}b_{i+1}=\prod_{s=i}^{k-1}\nfrac{(s-\mu_1)(s+1)}{(s-k)(n-s)}b_{k}; \\
&a_i=\textstyle\nfrac{i+1}{n-i}\prod_{s=i+1}^{k-1}\nfrac{(s-\mu_1)(s+1)}{(s-k)(n-s)}b_{k}, \qquad\qquad\qquad \text{where}\quad i=0,\dots, k-1;\\
&a_{k}=a_{k+1}=\cdots=a_{n-1}=(\mu_1-n)a_n=0;\\
&b_{k+1}=b_{k+2}=\cdots=b_{n-1}=0,
\end{split}
\ee
which depends on two parameters $(b_{k}, a_{n})$ if $\mu_1=n$ and on one parameter $b_{k}$ otherwise; hence, $v_{2n+1}$ either depends on $\nfrac{b_{k}}{a_{n}}\in\Cee\Pee^1$ or is defined uniquely (up to a~factor), respectively.

\underline{Case (ii)3. Let $\mu_2-n+k=0$ for some $k\in\{0,\dots, n-1\}$ and $\mu_1=t\in \{0,\dots, k-1\}$.} Thanks to the fact that $t<k$, we see that
$\mu_2-n+t\not=0$. Hence, $b_t=0$ and similarly, as in two previous cases, we obtain 
\be\label{Case3}
\begin{split}
a_{0}&=\cdots=a_{t-1}=0; \qquad \qquad b_{0}=\cdots=b_{t}=0;\\
a_{k}&=\cdots=a_{n-1}=a_n=0;\quad b_{k+1}=\cdots=b_{n-1}=0;\\
b_i&=\textstyle\nfrac{(i-\mu_1)(i+1)}{(\mu_2-n+i)(n-i)}b_{i+1}=\prod_{s=i}^{k-1}\nfrac{(s-\mu_1)(s+1)}{(s-k)(n-s)}b_{k}; \\
a_i&=\textstyle\nfrac{i+1}{n-i}\prod_{s=i+1}^{k-1}\nfrac{(s-\mu_1)(s+1)}{(s-k)(n-s)}b_{k}, \qquad\qquad \text{where}\quad i=t+1,\dots, k-1\\
\end{split}
\ee
which depends on two parameters $(a_{t}, b_{k})$; so, up to a~ factor, $v_{2n+1}$ depends on $\nfrac{a_{t}}{b_{k}}\in\Cee\Pee^1$.

\sssec{The answer} Recall that the $\fk(1|1)$-invariant differential operators between spaces of weighted densities are polynomials in $D_\theta:=\theta\del_t-\del_{\theta}$ with constant coefficients, not in $K_\theta=\theta\del_t+\del_{\theta}$, see \cite{LSh, Shch} and \cite[Theorem 3.2]{BLS1}.

\parbegin{Theorem}\label{ThGRC2} The classification of GRC operators in terms of $\fosp(1|2)$-singular vectors $v_k$ is given for the negatives of $\mu_1$ and $\mu_2$ by the formulas \eqref{Hek00}--\eqref{FinalEquations_01} for $k$ even and \eqref{Case1}--\eqref{Case3} for $k$ odd; mind that $(D_\theta)^2=-\partial_t$ and $f^{(i)}:=(\partial_t)^i(f)$:
\[
\begin{split}
\llb f,g\rrb_{2n}&=\textstyle\mathop{\sum}\limits_{0\leq i\leq n}\nfrac{1}{i!(n-i)!}c_if^{(i)}g^{(n-i)}\\
&\quad -\textstyle\mathop{\sum}\limits_{0\leq i\leq n-1}\nfrac{1}{i!(n-1-i)!}e_iD_\theta(f^{(i)}) D_\theta(g^{(n-i-1)});\\ 
\llb f,g\rrb_{2n+1}&=\textstyle\mathop{\sum}\limits_{0\leq i\leq n}\left(a_iD_\theta(f^{(i)})g^{(n-i)}+b_if^{(i)} D_\theta(g^{(n-i)})\right). 
\end{split}
\] \end{Theorem}

\parbegin{Example}\label{Ex1} We immediately see an analog of the following operator on the 1-dimensional manifolds, see \cite{Gr}, where $f':=\nfrac{\partial f}{\partial x}$,
\be\label{10}
\left(f, g\right)\longmapsto \llb f, g\rrb_1dx,\qquad \text{where \quad $\llb f, g\rrb_1:=
af'g+bfg'$}\qquad \text{for any $\nfrac{a}{b}\in\mathbb{CP}^1$ and $f,g\in\cF$}.
\ee
Indeed, thanks  to the fact \eqref{volISdtheta}, we can define
\[
[df]:=df\pmod{\cF\alpha}\equiv \left(\nfrac{\partial f}{\partial \theta}-\theta\nfrac{\partial f}{\partial t}\right) d\theta=D_\theta(f)d\theta,
\]
where $df:=dt\nfrac{\partial f}{\partial t}+d\theta\nfrac{\partial f}{\partial \theta}$. The operator
\be\label{11}
\begin{split}
\left(f, g\right)\longmapsto
&\left(a[df]g+bf[dg]\right)\\
=&\left(aD_\theta(f)g+(-1)^{p(f)}bfD_\theta(g)\right)d\theta\qquad \text{for any $\nfrac{a}{b}\in\mathbb{CP}^1$ and $f,g\in\cF$}
\end{split}
\ee
is invariant under the whole $\fk(1|1)$, as was noted in \cite{L1, GO}, hence is a~GRC operator.
\end{Example}

\ssec{The GRC operators and the algebra structure on the space of all weighted densities} Let $\Nee-\nfrac12:=\left\{\nfrac12, \nfrac32, \nfrac52, \dots\right\}$. Recall that $\Pi$ is the inversion of parity functor.

\sssec{Open problem} Set $\cF_{\bcdot}:=\oplus_w\ \cF_w$; find the coefficients $r_n$ and $s_n$ in the expressions
\be\label{A11}
f\ast g := \arraycolsep=1.5pt\begin{cases}\sum_{n\in\Nee}
r_n\llb f,g\rrb_n,  &\text{for any
$f, g\in \cF_{\bcdot}$};\\
\sum_{n\in\Nee-\nfrac12}
s_n\llb f,g\rrb_{2n},  &\text{for any $f, g\in \Pi(\cF_{\bcdot})$},\\
\end{cases}
\ee
to define associative multiplications on the spaces $\cF_{\bcdot}$ and~ $\Pi(\cF_{\bcdot})$. Conjecturally, $r_n=s_n=1$ for all $n$ is a~particular solution.

\section{The GRC operators on the general $(1|1)$-dimensional superstring} \label{SS11}

\ssec{On $\fvect(1|1)$-modules induced from $\fgl(1|1)$-irreducibles} \label{ssInd11}
Recall that a~\textit{superdomain} $\cU$ is a~pair $(U, \cF(U,V))=$(a~domain $U$, the superalgebra $\cF$ of functions on $U$ with values in the Grassmann algebra $\Lambda(V)$ of a~ vector space $V$); by definition, the \textit{morphisms of superdomains}$(U, \cF(U,V))\to (\widetilde U, \widetilde\cF:=\cF(\widetilde U,\widetilde V))$ are the parity preserving morphisms of the superalgebras $\cF\to \widetilde\cF$, see \cite{Lo} containing further details. The (super)\textit{dimension} of $\cU$ is the pair $(\dim U|\dim V)$.

Let the coordinates on the $(1|1)$-dimensional superdomain $\cU$ be $x$ (even) and $\xi$ (odd). Consider the \textit{standard $\Zee$-grading} of $\fg:=\fvect(1|1)$, namely, we set $\deg x=\deg \xi=1$. Hence, $\fg_0\simeq\fgl(1|1)$. Let $\del:=\del_x$ and $\delta:=\del_\xi$. For a~basis of $\fg_{0}=\fgl(1|1)$ realized as a~subalgebra of $\fvect(1|1)$, we take
\[
X_-:=\xi\del, \quad H_1:=x\del, \quad H_2:=\xi\delta,\quad \text{and}\quad X_+:=x\delta.
\]
The eigenvector with respect to $H_1$ and $H_2$ is called a~\textit{weight} vector, the weight vector annihilated by $X_+$ (resp., $X_-$) is called \textit{highest} (resp., \textit{lowest}). Let $V$ be a~$\fgl(1|1)$-module; setting $\fg_iV=0$ for $i>0$ we turn $V$ into a~module over $\fg_{\geq0}:=\oplus_{i\geq0}\ \fg_i$. The space of tensor fields of type $V$ is defined to be
\[
T(V):=\Hom_{U(\fg_{\geq0})}(U(\fg), V)\simeq\Cee[[x, \xi]]\otimes V.
\]
Hereafter, we work in terms of the space  dual to $T(V)$, i.e., we consider
\[
I(V^*):=U(\fg)\otimes_{U(\fg_{\geq0})} V^*\simeq\Cee[\del, \delta]\otimes V^*.
\]

Recall again that a \textit{singular vector} in ${I(V_1)\otimes I(V_2)}$  or $I(V)$ is a~weight vector (relative $H_1$ and $H_2$)   annihilated by ${\fg_{>0}:=\oplus_{i>0}\ \fg_i}$. A~ generalisation of Veblen's bilinear problem reduces to the task of classifying highest weight singular vectors in $I(V_1)\otimes I(V_2)$, where $V_1$ and $V_2$ are irreducible $\fgl(1|1)$-modules; this problem is solved, see a~review \cite{GLS2}. To describe $\fpgl(2|1)$-\textit{singular} vectors, see Subsection~\ref{B11}, is more cumbrous problem. We solve it in the subcase most resembling that of weighted densities, see Subsection~\ref{4cases}.
The other cases might be wild.



\ssec{A description of irreducible $\fgl(1|1)$-modules with highest (or lowest) weight vector}
 Let $V:=W^{a;b}$ be an irreducible module over the commutative Lie algebra $\fh=\Span (H_1, H_2)$, and $v\in V$ a~non-zero element such that $H_1v=av$ and $H_2v=bv$ for some $a,b\in\Cee$. Since $\fh$ is commutative, $\dim V=1$. Let $\fb:=\fh\oplus \Cee X_+$. For definiteness, let $v$ be \textit{even}. We make $V$ into a~$\fb$-module by setting $X_+V=0$ and set
\[
\begin{split}
I_{\fb}^{\fg_{0}}(V):&=U(\fg_{0})\otimes_{U(\fb)}V\simeq U(\Cee X_-\oplus \fb)\otimes_{U(\fb)}V\\
&\simeq U(\Cee X_-)\otimes U(\fb)\otimes_{U(\fb)}V\simeq\Cee [X_-]\otimes V
\\
&
{\simeq}(\Cee X_-\oplus\Cee 1)\otimes V\qquad\qquad\qquad\qquad\qquad\qquad (\text{since}\quad X_-^2=\textstyle\nfrac12[X_-, X_-]=0)\\
& \simeq(\Cee X_-\otimes V)\oplus V.
\end{split}
\]

Clearly, $I_{\fb}^{\fg_{0}}(V)$ is naturally graded by powers of $X_-$; let $\deg V:=0$. Let us find out when there exists a~highest weight vector of degree $-1$ annihilated by $X_+$:
\[
X_+X_-v= Ev= (a+b)v=0\quad \Rightarrow\quad a+b=0.
\]

It is convenient to change notation and set
\[
H:=H_1-H_2\quad \text{and}\quad E:=H_1+ H_2.
\]
 Accordingly, set
\be\label{ab}
V^{\lambda;\mu}:=W^{a;b} \qquad \text{for}\quad \lambda=a+b, \quad \text{and}\quad \mu=a-b.
\ee

Thus, if $\lambda\neq 0$, the module $M^{\lambda;\mu}:=I_{\fb}^{\fg_0}(V^{\lambda;\mu})$ is irreducible, while if $\lambda=0$, the module $M^{\lambda;\mu}$ is indecomposable, with a~submodule spanned by $X_-v$ and the quotient module spanned by $v$. The dual of this indecomposable module is also indecomposable with the submodule spanned by~ $v^*$. These indecomposable modules can be glued together in the same way as the adjoint $\fgl(1|1)$-module, whose maximal submodule is $\fsl(1|1):=\Span(X_-, E, X_+)$, is glued of 1-dimensional irreducibles. Note that $\fpgl(1|1):=\fgl(1|1)/\Cee E$. Let the actions of $X_-$ (resp., $X_+$) on the source of the arrow directed to the left (resp., right) and downwards be depicted as follows:
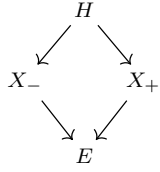
\[
\tiny
\begin{tikzcd}[column sep=tiny]
& H \ar[dl] \ar[dr] \\
X_- \ar[dr]&& X_+ \ar[dl] & \\
& E
&
&
\end{tikzcd}
\]
\normalsize

\subsection{The tensor product $\text{irr}(M^{\lambda;\, \mu})\otimes \text{irr}(M^{\sigma ; \, \rho })$} What are the irreducible $\fgl(1|1)$-modules ${M^{\lambda;\, \mu}\otimes M^{\sigma ; \, \rho}}$ is glued of? Recall that, by definition, the highest weight vectors $v\in M^{\lambda;\, \mu}$ and $w\in M^{\sigma ;\, \rho }$ are \textit{even}; if $v$ is odd, we denote the module it generates by $\Pi M^{\lambda;\, \mu}$. Let $\text{irr}(M^{\lambda;\,\mu})$ be the irreducible quotient of $M^{\lambda;\,\mu}$ modulo the maximal submodule; $\text{irr}(M^{\lambda;\,\mu})$ and $M^{\lambda;\, \mu}$ share the highest weight vector.

\sssec{The 4 cases}\label{4cases} Only the following cases are possible:

$(i)$ $\lambda=\sigma=0$. Clearly, $\text{irr}(M^{0;\, \mu})=V^{0;\, \mu}$ and $\text{irr}(M^{0; \, \rho })=V^{0;\, \rho }$; so $\dim V^{0;\, \mu}\otimes V^{0; \, \rho }=1$.

$(ii)$ $\lambda=0$ and $\sigma\neq 0$, then $V^{0;\,\mu}\otimes M^{\sigma; \, \rho }\simeq M^{\sigma;\, \rho +\mu}$ is irreducible. The case $\lambda\neq 0$ and $\sigma=0$ is analogous.

$(iii)$ $\lambda \sigma\neq 0$ and $\lambda +\sigma\neq 0$, then $M^{\lambda;\, \mu}\otimes M^{\sigma; \, \rho }$ is a~ direct sum of two irreducibles. Indeed, the vector
\be\label{subsp}
xX_- v\otimes w+yv\otimes X_- w\qquad \text{for some $x,y\in\Cee$}
\ee
is highest if
\[
xEv\otimes w+yv\otimes Ew=(x\lambda+y\sigma )v\otimes w=0.
\]
The last condition is true only for
$x=-y\nfrac{\sigma }{\lambda}$. Moreover,
\be\label{dir2}
M^{\lambda;\,\mu}\otimes M^{\sigma ; \, \rho }\simeq M^{\lambda+\sigma;\, \mu+\rho}\oplus \Pi M^{\lambda+\sigma;\, \mu+\rho-2}.
\ee

$(iv)$ $\lambda \sigma\neq 0$ and $\lambda +\sigma=0$, then $M^{\lambda;\, \mu}\otimes M^{\sigma ; \, \rho }$ is the sum of 4 modules (two even and two odd), which are trivial as $\fsl(1|1)$-modules but non-trivial as $\fgl(1|1)$-modules, glued in the same way as the adjoint $\fgl(1|1)$-module is glued of irreducibles, see diagram \eqref{simi}.

Comparing with the structure of submodules in the adjoint $\fgl(1|1)$-module, we'd expect that if $\lambda +\sigma =0$, there are TWO highest weight vectors on ``level $-1$", not one.
Additionally, we'd expect that the vector $X_- v\otimes X_- w$ is a~highest weight vector, whereas (recall that $X_+$ and $X_-$ are odd, and the Sign Rule)
\[
\begin{split}
X_+(X_- v\otimes X_- w)&= Ev\otimes X_- w-X_- v\otimes Ew\\
&=
\lambda v\otimes X_- w-\sigma X_- v\otimes w= \lambda (v\otimes X_- w+X_- v\otimes w)\neq 0.
\end{split}
\]

However, there is no mistake. 
Indeed, if $\lambda +\sigma =0$, then there is a~1-dimensional invariant subspace
spanned by the vector $v\otimes X_- w+X_- v\otimes w$: both $X_-$ and $X_+$ annihilate it. The quotient modulo this subspace has 2 invariant 1-dimensional subspaces spanned by
$v\otimes w$ and $X_-v\otimes X_-w$, see formulas:
\[
\begin{array}{lll}
X_-:&v\otimes w&\mapsto X_-v\otimes w+v\otimes X_-w;\\
X_+:& X_-v\otimes w-v\otimes X_-w&\mapsto 2\lambda(v\otimes w);\\
X_-:&X_-v\otimes w-v\otimes X_-w&\mapsto -2X_-v\otimes X_-w;\\
X_+: &X_-v\otimes X_-w&\mapsto \lambda(v\otimes X_-w+X_-v\otimes w).\\
\end{array}
\]
The quotient modulo the direct sum of these subspaces is also
1-dimensional, so the structure of submodules is similar to that of the adjoint $\fgl(1|1)$-module (indicated are the basis vectors and their weights relative $H$ on the right of the vector, the label at the arrow is the coefficient, if distinct from $1$, the target is multiplied by):
\be\label{simi}
\begin{tikzcd}[column sep=tiny]
& {\scriptstyle X_-v\otimes w-v\otimes X_-w\ (\mu+\rho-2)} \ar[dl, "-2"] \ar[dr, "2\lambda"] \\
{\scriptstyle X_-v\otimes X_-w\ (\mu+\rho-4)}\ar[dr, "\lambda"]&&{\scriptstyle v\otimes w\ (\mu+\rho)} \ar[dl] & \\
& {\scriptstyle X_-v\otimes w+v\otimes X_-w\ (\mu+\rho-2)}
&
&
\end{tikzcd}
\ee

\ssec{Conclusion}\label{wild} The only feasible case to tackle in the search of $\fpgl(2|1)$-singular vectors \eqref{Singfn}, \eqref{v_n} is the simplest Case $(i)$, see Subsection~\ref{4cases}, the classical setting of the GCR-problem. 

In Case $(i)$, we consider the tensor product of the $\fvect(1|1)$-modules $I^{0;\, \mu}$ induced from the 1-dimensional $\fgl(1|1)$-modules of the form $V^{0;\, \mu}$. Recall that the module of weighted densities $\cF_a=\cF\vvol^a$ over the algebra of functions $\cF$ on the superstring $\cC^{1|1}$ is spanned by an element $\vvol^a$, where $a\in\Cee$, on which the Lie derivative $L_D$ along the vector field $D$ acts as multiplication by $a\Div(D)$, where the divergence is defined by the formula
\be\label{Div}
\textstyle
\Div(f\del+g\delta):=\nfrac{\del f}{\del x}+(-1)^{p(g)}\textstyle\nfrac{\del g}{\del {\xi}}\qquad\text{~~for any $f,g\in\cF$}.
\ee

The  dualization $(\cF_{-\mu})^*\simeq I^{0;\,\mu}$ gives an interpretation of the results of Subsection~\ref{Sres}.
\ssec{Problem B on a~$1|1$-dimensional superdomain}\label{B11} Consider ${\fpgl(2|1)\subset \fvect(1|1)}$ embedded as a~graded subalgebra so that $\fpgl(2|1)_i=\fvect(1|1)_i$ for $i=-1$ and $0$ in the standard grading.

For a~basis of $\fvect(1|1)_1$ we consider two vectors $s_x:=xE$ and $s_\xi:=\xi E=x\xi\del$, where $E=x\del+\xi\delta$ belongs to $\fpgl(2|1)\subset \fvect(1|1)$, and two divergence-free elements which do not belong to $\fpgl(2|1)$, so we do not need them in this note.

Let $V_1$ and $V_2$ be 1-dimensional $\fg_0$-modules. Set $I(V_1):=\Cee[\del',\delta']\otimes V_1$ and $I(V_2):=\Cee[\del'',\delta'']\otimes V_2$, where $\del'$ and $\delta'$ (resp., $\del''$ and $\delta''$) are copies of $\del$ and $\delta$, as Grozman suggested in \cite{Gr}, and use the fact that
\[
I(V_1)\otimes I(V_2)\simeq \Cee[\del',\delta',\del'',\delta'']\otimes V_1\otimes V_2.
\]
Then,  any weight (with respect to $H_1$ and $H_2$) vector of $I(V_1)\otimes I(V_2)$ on level $-n$ is of the form
\be\label{Singfn}
\begin{split}
v_n=\sum_{0\leq k\leq n} &(\del')^k (\del'')^{n-k} a_k \otimes r_k+\sum_{0\leq k\leq n-1} (\del')^{k} (\del'')^{n-k-1}\delta' \ b_k\otimes s_k \\[2mm]
&+
\sum_{0\leq k\leq n-1} (\del')^k(\del'')^{n-k-1} \delta'' \ c_k\otimes t_k+\sum_{0\leq k\leq n-2} (\del')^{k} (\del'')^{n-k-2} \delta' \delta'' \ d_k\otimes u_k,
\end{split}
\ee
where $ a_k, b_k, c_k, d_k\in V_1$ and $r_k, s_k, t_k, u_k\in V_2$.

Note that $\fg_>$ is generated by $s_{\xi}$ and one more (divergence-free) element. As we showed above, the description of Gordan-Rankin-Cohen operators is equivalent to the description of the $\fpgl(2|1)$-\textit{singular} vectors, i.e., the highest weight vectors killed only by $s_{\xi}$. 
 Thus, our task is to solve the system
\be\label{sys}
\begin{cases}
X_+(v_n)=0; \\
s_\xi(v_n)=0.
\end{cases}
\ee
Moreover, it suffices to confine ourselves to vectors
\be\label{v_n}
v_n\in \Cee[\del',\delta',\del'',\delta'']\otimes \text{irr}(V_1\otimes V_2),
\ee
where $\text{irr}(V_1\otimes V_2)$ is an irreducible $\fgl(1|1)$-component (submodule) of $V_1\otimes V_2$, where the $V_i$ are irreducible $\fgl(1|1)$-modules. Since $\dim(V_1)=\dim(V_2)=1$, it follows that $\text{irr}(V_1\otimes V_2)=V_1\otimes V_2$.


\ssec{The $\fpgl(2|1)$-singular vectors in Case $(i)$}
\label{Sres}
Let $H_1'$ (resp., $H_1''$) designate the operator $H_1$ acting on the 1st (resp., 2nd) factor in the tensor product $u:= v_1\otimes v_2\in V_1\otimes V_2$, where $V_1:=V^{0;\mu_1}$ and $V_2:= V^{0; \mu_2}$. Thus, we seek $\fpgl(2|1)$-singular vectors of the form
\be\label{Singfn_2}
\begin{array}{ll}
v_n=\Big(\sum_{k=0}^{n}& A_k(\del')^k (\del'')^{n-k} +\sum_{k=0}^{n-1} B_k(\del')^{k} (\del'')^{n-k-1}\delta' \\[2mm]
&+
\sum_{k=0}^{n-1} C_k(\del')^k(\del'')^{n-k-1} \delta'' +\sum_{k=0}^{n-2} D_k(\del')^{k} (\del'')^{n-k-2} \delta' \delta'' \Big)u.
\end{array}
\ee

To simplify the task, we use two observations due to A.~Lebedev. First, the summands of~ $v_n$ can be divided into the even (with coefficients $A$ and $D$) and odd (with coefficients $B$ and $C$), to be considered separately. Second, the even summands can be further subdivided: since all summands in $v_n$ are homogeneous under the action of $H_2$, then the spaces of singular vectors are also direct sums of subspaces homogeneous with respect to parity and $H_2$-weight. Hence, we can consider the summands with coefficients~ $A$ (of $H_2$-weight~ 0) and the summands with coefficients~ $D$ (of $H_2$-weight~ $-2$) separately.
Thus, the space of vectors of degree $-n$ can be subdivided into $3$ types of subspaces:
\begin{itemize}

\item[(A)] with basis vectors $(\del')^k (\del'')^{n-k} u$, where $k=0,\dots, n$;

\item[(D)] with basis vectors $(\del')^k (\del'')^{n-k-2}\delta'\delta'' u$, where $k=0,\dots, n-2$;

\item[(BC)] with basis vectors $(\del')^k (\del'')^{n-k-1}\delta' u$, and $ (\del')^k (\del'')^{n-k-1}\delta'' u$, where $k=0,\dots, n-1$.
\end{itemize}

\underline{In case (A)}, we have
\[
X_+ (\del')^k (\del'')^{n-k} u = -k (\del')^{k-1} (\del'')^{n-k}\delta' u - (n-k)(\del')^k (\del'')^{n-k-1}\delta'' u.
\]
The summands of images of different basis elements are also different basis elements, and none of the images vanishes, except for the case of $n=0$. Hence, the kernel of $X_+$ on the space of $A$-vectors is zero.

\underline{In case (D)}, we have
\be\label{CaseD}
\begin{split}
\textstyle
X_+ (\del')^k (\del'')^{n-k-2}\delta'\delta'' u =& 0;\\
s_\xi (\del')^k (\del'')^{n-k-2}\delta'\delta'' u =& (-k + \nfrac12\mu_1)(\del')^k (\del'')^{n-k-2}\delta'' u\\
& - (-(n-k-2) + \textstyle\nfrac12\mu_2)(\del')^k (\del'')^{n-k-2}\delta' u.
\end{split}
\ee
The images of different basis elements are still different basis elements, but in this case, the image of the $k$-th basis element can vanish if $-k +\nfrac12 \mu_1$ and $-(n-k-2) + \nfrac12\mu_2$ vanish simultaneously.
Such $k$ exists if and only if $\mu_1$ and $\mu_2$ are non-negative even numbers and $\mu_1+\mu_2=2n-4$.
We have
\begin{equation}\label{Ev}
\text{\begin{minipage}[c]{14cm} If $\mu_1=2k$ for some $k=0, 1, \dots, n-2$, then $\mu_2=2n-4-2k$ and $D_k\in\Cee$ is arbitrary whereas $D_i=0$ for $i\neq k$.\end{minipage}}
\end{equation}

\underline{In case (BC)}, we have
\[
\begin{split}
X_+ (\del')^k (\del'')^{n-k-1}\delta' u &= (n-k-1)(\del')^k (\del'')^{n-k-2}\delta'\delta'' u; \\
X_+ (\del')^k (\del'')^{n-k-1}\delta'' u &= -k(\del')^{k-1} (\del'')^{n-k-1}\delta'\delta'' u.
\end{split}
\]
This implies that the dimension of the kernel of $X_+$ on the space of BC-vectors is equal to $(0|n+1)$; and for a~basis of this space we can take
\be\label{els}
k(\del')^{k-1} (\del'')^{n-k}\delta' u + (n-k)(\del')^k (\del'')^{n-k-1}\delta'' u, \qquad \text{where $k=0,\dots, n$.}
\ee
Note that certain elements of this basis involve $(\del')^{-1}$ and $(\del'')^{-1}$ which appear only with zero coefficients, so these expressions make sense.

Now, let us see how $s_\xi$ acts on these elements \eqref{els}:
\[
\begin{split}
s_\xi (k(\del')^{k-1} (\del'')^{n-k}\delta' u &+ (n-k)(\del')^k (\del'')^{n-k-1}\delta'' u) \\
=& k(-(k-1)+\textstyle\nfrac12\mu_1)(\del')^{k-1} (\del'')^{n-k} u \\
&+ (n-k)(-(n-k-1)+\textstyle\nfrac12\mu_2)(\del')^k (\del'')^{n-k-1} u.
\end{split}
\]
In other words, if an element of the kernel of $X_+$ on the BC-space is of the form
\[
b=\sum_{0\leq k\leq n} b_k\left(k(\del')^{k-1} (\del'')^{n-k}\delta' u + (n-k)(\del')^k (\del'')^{n-k-1}\delta'' u\right),
\]
then the condition $s_\xi b = 0$ boils down to the conditions we obtain looking at the coefficient of $(\del')^k (\del'')^{n-k-1} u$ in $s_\xi b$: 
\be\label{eqs}
\begin{split}
(n-k)(\mu_2-2&(n-k-1)) b_k \\&+ (k+1)(\mu_1-2k) b_{k+1} = 0, \qquad \text{where}\quad
k=0,\dots, n-1.
\end{split}
\ee

We solve this system in the same way as we solved system \eqref{Sect2}.
Denote by $M$ the $n\times (n+1)$ matrix whose elements are the coefficients of system \eqref{eqs}.
Clearly, $M$ can have non-zero elements only on the main diagonal
$D=(A_0,\ldots, A_{n-1})$, where $A_i=(n-i)(\mu_2-2(n-i-1))$, and on the next diagonal just above the main one, $(B_0,\ldots, B_{n-1})$, where $B_i=(i+1)(\mu_1-2i)$. Obviously, at most one $A_i$ (resp., $B_i$) can vanish. To solve a~system $A_i b_i + B_i b_{i+1} = 0$, where $ i=0, \ldots , n-1$ for $n+1$ unknowns $b_0, \ldots , b_n$, where at most one of the coefficient $A_i$ vanishes and at most one of the coefficient $B_i$ vanishes, it suffices to consider the 2 cases:
 
Case 1: $A_i = 0$ for some $i$, and $B_j = 0$ for some $j$ such that $j\geq i$.
In this case, the space of solutions is 2-dimensional; it is spanned by the following two solutions:

Solution 1: 
\be\label{sol1}
\text{$b_m = (-1)^m \prod_{k=0}^{m-1}A_k \prod_{k=m}^{i-1}B_k$ \qquad  for $m=0, \ldots , i$; 
and $b_m = 0$ for $m = i+1, \ldots , n$.}
\ee

Solution 2: 
\be\label{sol2}
\text{$b_m = 0$ for $m = 0, \ldots , j$; and $b_m = (-1)^m \prod_{k=j+1}^{m-1} A_k \prod_{k=m}^{n-1}B_k$ for $m = j+1, \ldots , n$.
}
\ee

Case 2: all other possibilities (either none of the $A_i$ and $B_j$ vanish; or some $A_i$ vanishes but none of the $B_j$ vanishes; or some $B_i$ vanishes but none of the $A_i$ vanishes; or some $A_i$ vanishes and some $B_j$ vanishes, where $j<i$.)
In this case, the space of solutions is 1-dimensional and is spanned by the following solution:
\be\label{solCase2}
b_m = (-1)^m \prod_{k=0}^{m-1}A_k \prod_{k=m}^{n-1}B_k \qquad \text{for $m=0, \ldots , n$.}
\ee
If some of the $A_i$ and $B_j$ vanish, then some of the $b_m$ in this solution also vanish, but not all (under the above conditions), i.e., this is a~non-trivial (non-zero) solution.

\sssbegin{Theorem}\label{ThMain} The dimension of the space of singular vectors $v_n$ is equal to
\begin{itemize}
\item[(i)] $1|1$ if $\mu_1$ and $\mu_2$ are non-negative even numbers and $\mu_1+\mu_2=2n-4$;
\item[(ii)] $0|2$ if $\mu_1$ and $\mu_2$ are even integers between $0$ and $2n-2$ inclusive, and $\mu_1+\mu_2\geq 2n-2$;
\item[(iii)] $0|1$ otherwise,
\end{itemize}
where, up to a~non-zero constant factor, the odd singular vectors $v_n$ are uniquely defined in case $(iii)$ and case $(i)$ by the formula~\eqref{solCase2}; in case $(ii)$, there is a~$1$-parameter family of odd $\fpgl(2|1)$-singular vectors given by formulas \eqref{sol1} and \eqref{sol2}; in case $(i)$, there is also an even $\fpgl(2|1)$-singular vector given by formula~\eqref{Ev}.
\end{Theorem}

\sssec{Open problem} Set $\cF_{\bcdot}:=\oplus_w\ \cF_w$; find the coefficients $r_n$ and $s_n$ in the expressions
\be\label{A111}
f\ast g := \arraycolsep=1.5pt\begin{cases}\sum_{n\in\Nee}
r_n\llb f,g\rrb_{g,n},  &\text{for any $f, g\in \cF_{\bcdot}$ and $\llb -,-\rrb_n$ even};\\
\sum_{n\in\Nee}
s_n\llb f,g\rrb_{g,n},  &\text{for any $f, g\in \Pi(\cF_{\bcdot})$ and $\llb -,-\rrb_n$ odd},\\
\end{cases}
\ee
to define associative multiplications on the spaces $\cF_{\bcdot}$ and~ $\Pi(\cF_{\bcdot})$.

\section{Discussion: on the two problems \eqref{2prob}}

\textbf{Solution to problem A}. The classical GRC \textit{brackets} appear in number theory, representation theory, in functional analysis, and in theoretical physics; for useful references, see \cite{Z, Z1, Z2, Gi, GTh, KP2, RTY, BSCK, BSCK1, Do} and especially \cite{KP1} containing the most complete results on bilinear differential invariant operators between the spaces of sections of homogeneous bundles over the six possible types of symmetric spaces. 

\textbf{Solution to problem B}. To classify GRC \textit{operators} one needs only simple linear algebra, as in this article, whereas to classify GRC \textit{brackets} requires a~more complicated technique, see \cite{Gi, CMZ, Z2}.

\ssec{Further open problems and conjectures} 1) In the case of 1-dimensional varieties over the ground field of characteristic $p>0$, all bilinear invariant operators between spaces of weighted densities are classified, see \cite{BL}, thus providing with examples of GRC operators as well. To list all GRC operators invariant under the subalgebra $\fsl(2)$ in the Lie algebra of vector fields seems to be a~ feasible open problem; 
however, unlike the open problems listed in \cite{BL}, no interpretation of this problem is known at the moment.

2) In an unfinished draft of the description of GRC operators between spaces of weighted densities on $(1|N)$-dimensional superstrings with a~ contact structure (by Bouarroudj et al.) many interesting new $\fk(1|N)$-invariant operators, hence GRC operators, are found. In particular, several incomplete for $N=2$ results of \cite{BLO} dealing with problem B are completed: (i)~ for $N=2$, the weight is a~\textit{pair} of numbers, (ii)~ a~contact structure with odd time is possible --- both cases were never previously considered. We hope Bouarroudj will finish this draft.

3)  In \cite{OR}, analogs of GRC operators invariant under $\fo(n+2)\subset\fvect(n|0)$ between spaces of weighted densities are described (not classified) for a~real form of $\fo(n+2)$; for the classification of GRC operators invariant under $\fsl(2)\simeq\fo(3)\subset\fvect(1|0)$ over $\Cee$, see \cite{BoLe}. 

It seems also feasible to classify GRC operators invariant under other maximal finite-dimensional simple Lie sub(super)algebras of simple infinite-dimensional vectorial Lie (super)algebras.  However, if we generalize the setting of problem B by considering binary operators invariant under $\fpgl(a+1|b)\subset \fvect(a|b)$ between spaces of tensor fields with multidimensional fibers, then Grozman's classification of the $\fvect(a)$-invariant operators, see \cite{Gr}, gives us a~part of the answer. To tackle the problem in such generality is hardly feasible, see Grozman's proof and Subsection~\ref{wild}.

\subsection*{Disclosures}
We are thankful to  I.~Shchepochkina  and A.~Lebedev for help.
The  research of V.~Bovdi was supported by UAEU grant G00004159. 
No conflict of interest.
The data used to support the findings of this study are included within the article.


\def\eightit{\it}
\def\bib{\bf}
\bibliographystyle{amsalpha}

\end{document}